\documentclass[reqno]{amsart}
\usepackage[marginratio=1:1]{geometry}

\usepackage[hidelinks]{hyperref}
\usepackage{calc}
\newsavebox\CBox
\newcommand\hcancel[2][0.5pt]{%
  \ifmmode\sbox\CBox{$#2$}\else\sbox\CBox{#2}\fi%
  \makebox[0pt][l]{\usebox\CBox}%
  \rule[0.5\ht\CBox-#1/2]{\wd\CBox}{#1}}
\usepackage{blindtext}
\usepackage{color}
\usepackage{hyperref,url}
\usepackage{comment}
\usepackage{float}
\usepackage{hyperref}
\usepackage{enumerate}
\usepackage{enumitem}   
\usepackage{bbm}
\usepackage{cancel}
\usepackage{mathrsfs}
\usepackage{colonequals}
\usepackage{pdfpages}

\usepackage{amsmath,amsfonts,amsthm,bm}
\usepackage{mathrsfs}  
\usepackage{xcolor}
\usepackage{amsfonts}
\usepackage{amssymb}
\usepackage{amsmath}
\usepackage{mathtools}
\usepackage{mathrsfs}  

\usepackage[normalem]{ulem}

\numberwithin{equation}{section}
\usepackage[toc,page]{appendix}
\usepackage{tikz-cd}
\theoremstyle{definition}

\newtheorem{definicao}{Definition}[section]
\newtheorem{remark}[definicao]{Remark}

\theoremstyle{plain}

\newtheorem{teorema}[definicao]{Theorem}
\newtheorem{suposicao}{Hypothesis}

\newtheorem{lema}[definicao]{Lemma}
\newtheorem{proposicao}[definicao]{Proposition}
\newtheorem{corolario}[definicao]{Corollary}

\usepackage{scalerel}
\newenvironment{TheoremproofD}[1]{\par\noindent{\emph{Proof of Lemma \ref{simp}}} \space#1}{\leavevmode\unskip\penalty9999 \hbox{}\nobreak\hfill\quad\hbox{$\qed$}}

\definecolor{roxo}{rgb}{0.44, 0.16, 0.39}
\definecolor{ao(english)}{rgb}{0.0, 0.5, 0.0}
\definecolor{dmagenta}{RGB}{139, 0, 139}
\definecolor{dgreen}{RGB}{0,90,0}
\definecolor{navy}{RGB}{0,0,128}

\usepackage{stackengine}

\def\d{\mathrm d}
\def \d {\mathrm{d}}
\def \supp {\mathrm{supp}}

\definecolor{iblue}{RGB}{0, 35, 194}
\title[Existence of horseshoes]{Nonuniform expansion and diffusive noise imply random horseshoes and random Young towers} 

\address{$^{2}$International Research Center for Neurointelligence, The University of Tokyo, Tokyo 113-0033, Japan}
\email{\href{mailto:jeroen.lamb@imperial.ac.uk}{jeroen.lamb@imperial.ac.uk}}
\author[Jeroen S.W.~Lamb]{Jeroen S.W. Lamb$^{1,2}$}
\address{$^{1}$Department of Mathematics, Imperial College London, London SW7 2AZ, UK}
\author[Giuseppe Tenaglia]{Giuseppe Tenaglia*$^{1}$}
\email{\href{mailto:giuseppe.tenaglia20@imperial.ac.uk}{giuseppe.tenaglia20@imperial.ac.uk}}
\author[Dmitry Turaev]{Dmitry Turaev$^{1}$}
\email{\href{mailto:d.turaev@imperial.ac.uk}{d.turaev@imperial.ac.uk}}

\allowdisplaybreaks
\begin{document}



\begin{abstract}
{We propose  a notion of random horseshoe and prove density of random horseshoes 
and existence of a random Young tower with annealed exponential tails for nonuniformly expanding random dynamical systems with diffusive noise.}
\end{abstract}
\maketitle
   
\section{Introduction}
\label{intro}
A random dynamical system is a dynamical system  whose time evolution depends both on the current state and on an underlying source of randomness. Like deterministic dynamical systems, these may display chaos in the sense of sensitive dependence of trajectories on initial conditions.

The study of chaotic random dynamical systems has developed along two main perspectives. The first is a statistical approach focussing on stationary measures and the positivity of their associated Lyapunov exponents \cite{blumenthal2022positive,blumenthal2017lyapunov,Lian2012PositiveLE}. 
The second perspective is 
topological, focussing on 
the relative behaviour of trajectories with different initial conditions for 
typical noise realisations, in particular in the case when the top Lyapunov exponent is positive.

For the class of nonuniformly expanding  deterministic chaotic systems, a well-developed topological theory is available. These are systems admitting an ergodic invariant measure, absolutely continuous with respect to volume, whose Lyapunov exponents are positive. 
Such systems feature fundamental topological structures, such as horseshoes and Young towers,
which are widely regarded as hallmarks of chaotic dynamics. 

For  an  endomorphism  $f \colon M \to M$,  two disjoint boxes $B_0,B_1 \subset M$ form a horseshoe if there exists $N \in \mathbb{N}$ such that
\begin{align*}
f^N(B_0) \cap f^{N}(B_1) \supset B_0 \cup B_1,
\end{align*}
and $f^N$ is uniformly expanding
on the preimage-set 
\begin{align*}
S:=\{x \in B_0 \cup B_1 \colon f^N(x) \in B_0 \cup B_1 \},
\end{align*}
i.e. there exists $\lambda >1$ such that for all $x\in S$
\begin{align*}
|df^N(x)\cdot v| \ge \lambda |v| \qquad \forall v \in T_xM.
\end{align*}
This implies that on the maximal invariant set 
\begin{align*}
H = \bigcap_{j \in \mathbb{N}} f^{-jN}(B_0 \cup B_1)
\end{align*}
  $f^N\mid_{H}$ is uniformly expanding and topologically conjugate to the one-sided shift on two symbols, so that the infinite sequences of zeros and ones are in one-to-one correspondence with the orbits of points in $H$: 
\begin{equation}\label{simboli}
\forall  \{s_n\}_{n \in \mathbb{N} } \in \{0,1\}^{\mathbb{N}},\,\,\exists\,!\, x \in H \colon  f^{Nj}(x) \in B_{s_j} \qquad \forall j \in \mathbb{N}.
\end{equation}
The topological conjugacy implies that $f^N\mid_{H}$ has positive topological entropy and is sensitively dependent on initial conditions.
Similar constructions can be made by taking $k>2$ boxes.

Young towers  may be viewed as the countable-state analogy of horseshoes.  More specifically, let $m$ be the Lebesgue measure on $M$. Suppose there exists a reference set $\Delta_0 \subset M$ with $m(\Delta_0)>0$ satisfying the following properties:
\begin{itemize}
\item there exists a countable mod $0$  partition $\mathcal{P}$ of $\Delta_0$ into disjoint invertibility domains of f;
\item there exists a function $R \colon \Delta_0 \to \mathbb{N}$, such that $R$ is constant on every element of $\mathcal{P}$, and if $w \in \mathcal{P}$, then  $f^{R(w)}(w) =  \Delta_0\, \mod 0$.
\end{itemize}
Then, we define the induced map $F \colon \Delta_0 \to \Delta_0$ on each $w \in \mathcal{P}$ as
\begin{align*}
F_{|w} = f^{R(w)}_{|w}.
\end{align*}
If this map is uniformly expanding, then its discrete-time suspension with roof function $R$ is called a Young tower. These are extensively studied as they provide a framework for deriving a wide range of statistical properties, since the pioneering work of Lai-Sang Young \cite{Young1998,Young1999}.

In the deterministic setting, it is well known that non-uniformly expanding systems admit horseshoes and Young towers if trajectories do not return too frequently near the region of points where the derivative of $f$ vanishes, a condition commonly referred to as \emph{slow recurrence to the critical set} \cite{gelfert2010repellers,Alves}.

It is natural to ask whether the existence of analogous structures can be proved in the corresponding random setting. The key issue is that to build random horseshoes and random Young towers, one needs to establish  uniform geometric control of return times  across  typical realizations, based on the ergodic properties of the noise-averaged process.
This is a non-trivial  issue which simply does not appear in the deterministic settings.

In previous works \cite{lamb2023horsehoes,castro2023random}, we resolved the above issue for predominantly expanding one-dimensional random circle endomorphisms. These maps are strongly expanding throughout the whole phase space, except for a  comparatively small region in which the derivatives can approach zero. For this class, \cite{lamb2023horsehoes} shows the  density of random horseshoes, and  \cite{castro2023random} proves the  existence of a random Young-towers with annealed exponential tails.

 Random Young towers were constructed adapting the deterministic methods for small random perturbations of transitive non-uniformly expanding maps \cite{AST_2003__286__25_0,ALVES_VILARINHO_2013,su2023random}, and for certain specific families of random endomorphisms \cite{ASENS_2002_4_35_1_77_0,bahsoun2019quenched}. Random horseshoes were built for  mixing uniformly hyperbolic random dynamical systems \cite{huang2016ergodic, WenHUANG:281}, and 
 random non-hyperbolic  horseshoe-like constructions were developed in \cite{https://doi.org/10.1002/cpa.21698,huang2023full,huang2023observable}.



In this paper, under the sole assumption of diffusive noise, we establish density of random horseshoes for \(C^{1+\beta}\) multidimensional non-uniformly expanding random dynamical systems with non-degenerate critical points (Theorem \ref{thm1}). For this class, the estimates obtained to prove Theorem \ref{thm1} also allow us to carry out verbatim the tower construction in \cite{castro2023random}  and establish the existence of random Young towers  with annealed exponential tails (Theorem \ref{thm2}).

In this way, we establish a chain of arguments that connects, under minimal assumptions, random non-uniform expansion to  geometric features. Two key results underlie this. The first one is that diffusive noise implies  large-deviation estimates for noise-dependent unbounded observables  (see Section \ref{section3}). The second is that large-deviation estimates imply, combined with ergodicity, exponential tails for a certain sequence of return times with controlled geometry, the so-called Young times (see Section \ref{section5}). 

In the proof of Theorem \ref{thm1}, we show that the exponential tails of the  Young times are a sufficient condition for the density of random horseshoes. In Theorem \ref{thm2}, the exponential tails of Young times imply the existence of a random Young tower with annealed  exponential tails. 


From a technical point of view, it is essentially the diffusive noise assumption that allows us to remove the predominant-expansion assumption. Indeed, in \cite{lamb2023horsehoes,castro2023random}, large-deviation estimates were established using predominant 
expansion. Here, instead, we use the smoothing properties of the diffusive noise to exploit the perturbative approach developed in \cite{bougerol1985products}, and  prove,  in Proposition \ref{devmultid}, the central large-deviation result. 
The key challenge is that the observables are both noise-dependent and unbounded. To deal with these issues,  we have to consider an augmented Markov chain and carry out the analytic perturbation argument in a Banach space different from the one used in  \cite{bougerol1985products}.
 
In Section \ref{section6}, we develop a  probabilistic and combinatorial argument  to deduce the density of random horseshoes directly from the exponential tails of Young times.  This is a substantial improvement 
of the approach in \cite{lamb2023horsehoes}, where   the existence of a full-branch interval was essential.

This paper is organized as follows:  In Section \ref{Section2}, we recall the basic setup of random dynamical systems, state our hypotheses and main results. In Section \ref{section3} we establish large deviations results for the Lyapunov exponent and for recurrences to the critical set. 
In Section \ref{section5}, we review  the main properties of Young times and obtain the desired control on the statistics of Young times. In Section \ref{section6} we deduce, from the properties of Young times, the density of the random horseshoes. Finally, in Section \ref{section7}, we prove Theorem \ref{thm2}.

\section{Hypotheses and main results}
\label{Section2}
\subsection{Hypotheses}
Let $M$ be a smooth compact $k$-dimensional manifold  with Lebesgue (volume) measure $m$. Take $\sigma>0$ and consider the cube $[-\sigma,\sigma]^k$. Let $\Sigma$ be the Borel sigma-algebra on the cube and $\nu = \frac{\text{Leb}_{|[-\sigma,\sigma]^k}}{(2\sigma)^k}$. Consider the product space $\Omega:= \left([-\sigma,\sigma]^{k}\right)^{ \mathbb{N}}$, equipped with the product measure $\mathbb{P} = \left(\nu\right)^{\otimes \mathbb{N}}$ and the product $\sigma$-algebra $\mathcal{F}=\Sigma^{\otimes\mathbb{N}}$. The shift map
\begin{align*}
\theta \colon \Omega \to \Omega, \qquad 
\theta(\{\omega_i\}_{i \ge 0}) = \{\omega_{i+1}\}_{i \ge 0}
\end{align*}
is ergodic with respect to $\mathbb{P}$.

We consider a  dynamical system of the form
\begin{equation}\label{skewproduct}
\Theta \colon M \times \Omega \ \to   M \times \Omega , \qquad  
\Theta(x,\omega):= (f_{\omega}(x),\theta(\omega)),
\end{equation}
where $\{f_{\omega} \}_{\omega \in \Omega}$ is a family of maps such that
\begin{equation}\label{additiveNoise}
f_{\omega}(x) = g(x,\omega_0),
\end{equation} 
where $g \colon  M\times [-\sigma,\sigma]^k \to M$. Note that the random maps $f_{\omega}$ depend on $\omega$ only via  $\omega_0$.

We will use the following notation for the random iterations:
\begin{equation}\label{rdsiterations}
f^n_\omega(x)=f_{\theta^{n-1}\omega}\circ\ldots\circ f_\omega(x).
\end{equation}
\renewcommand\thesuposicao{(H1)}
Our first hypothesis concerns the regularity of $g$.
\begin{suposicao} \label{(A)} 
The function $g$ is  continuous and, for each $\omega_0\in[-\sigma,\sigma]^k$, $g(\cdot,\omega_0)$ is $C^{1+\beta}$, $\beta>0$, outside a critical
set $\mathcal{C}_{\omega_0}\subset M$ where this map is not differentiable or the derivative is not invertible.
Moreover,  for all $x \in M$, the function
$
g(x,\cdot) \colon [-\sigma,\sigma]^k \to M
$
is a diffeomorphism onto its image, and the map $(x,\omega_0) \to \det d_{\omega_0} g(x,\omega_0)$ is Holder-continuous. In particular, there exist $c_1,c_2>0$ such that
\begin{equation}\label{boundsforinvmeas}
c_1\ge \left|\det \left( d_{\omega_0} g(x,\omega_0)\right)\right|\ge c_2 >0\qquad \forall (x,\omega_0) \in M \times [-\sigma,\sigma]^k.
\end{equation}
\end{suposicao}
Hypothesis \ref{(A)} encompasses our notion of diffusive noise, and is crucial to establish large deviation estimates for the Lyapunov exponent (see Proposition \ref{devmultid}).
\begin{remark}\label{notation}
In view of \eqref{additiveNoise}, we use the notations \(g(x,\omega_0)\) and \(f_\omega(x)\) interchangeably. Similarly, we will write \(d_xg(x,\omega_0)\) and \(df_\omega(x)\) interchangeably for the derivative of the random maps with respect to \(x\).    
\end{remark}

Our second assumption concerns the regularity of the critical sets $\mathcal{C}_{\omega_0}$.

\renewcommand\thesuposicao{(H2)} 
\begin{suposicao} \label{(B)} 

For all $\omega_0 \in [-\sigma,\sigma]^k$,  the critical set $\mathcal{C}_{\omega_0}$ is a closed set with zero Lebesgue measure  and satisfies the following regularity conditions: there exist constants $B>1$ and $\alpha$ and $\beta>0$ as in Hypothesis \ref{(A)}  such that
\begin{enumerate}
\item  for all $x \in M$ and $\omega_0 \in [- \sigma,\sigma]^k$ 
\begin{equation}\label{c1}
\frac{1}{B} \text{dist}(x,\mathcal{C}_{\omega_0})^{\alpha}\le \frac{|d_xg(x,\omega_0)\cdot v|}{|v|} \le B \text{dist}(x,\mathcal{C}_{\omega_0})^{-\alpha},\qquad \forall  0\neq v \in T_x M
\end{equation} 
\item for all  $x,y \in M \setminus \mathcal{C}_{\omega_0}$ with  $|x-y| \le \frac{\text{dist}\left(x,\mathcal{C}_{\omega_0}\right)}{2}$
\begin{equation}\label{c2}
\left|\log||d_xg(x,\omega_0) ^{-1}||-\log||d_xg(y,\omega_0) ^{-1}|| \right| \le \frac{B}{\text{dist}\left(x,\mathcal{C}_{\omega_0}\right)^{\alpha}}|x-y|^{\beta}
\end{equation}
and 
\begin{equation}\label{c3}
\left|\log|\det(d_xg(x,\omega_0))|-\log|\det(d_xg(y,\omega_0))|\, \right| \le \frac{B}{\text{dist}\left(x,\mathcal{C}_{\omega_0}\right)^{\alpha}}|x-y|^{\beta}.
\end{equation}
\end{enumerate}
Furthermore, the map $\omega_0 \to \mathcal{C}_{\omega_0}$ is upper semi-continuous in the Hausdorff topology.
\end{suposicao}

If $\mu$ is any stationary measure for $f_{\omega}$, then it follows from equation \eqref{boundsforinvmeas} in Hypothesis \ref{(A)} and the estimates for $\mu$  in \cite{Lian2012PositiveLE} that  $\mu$ is absolutely continuous with respect to $m$ with uniformly bounded density: there exists a constant $C_0>0$ such that
\begin{equation}\label{BVestimate}
\left|\left|\frac{d\mu}{dm} \right|\right|_{\infty} \le C_0.
\end{equation}

Furthermore,  Hypothesis \ref{(A)} implies, by 
\cite[Theorem 3]{Lin1975} and \cite[Theorem B]{CastroLambOliconMendezRasmussen2024},
that there are at most a finite number of ergodic measures, and each of them carries at most finitely many cycling components.

Given $(x,\omega) \in M \times \Omega$, we define the Lyapunov exponent
\begin{align*}
\lambda(x,\omega):= -\limsup_{n \to \infty} \frac{1}{n} \sum_{i=0}^{n-1}\log\left|\left|\left(df_{\theta^i(\omega)}(f^i_{\omega}(x)\right)^{-1} \right| \right|.
\end{align*}

If $\mu$ is  ergodic   and
$\log|| df_{\omega}(x)^{-1}|| \in L^1(m\otimes \mathbb{P})$, then, because of\eqref{BVestimate},  $\log|| df_{\omega}(x)^{-1}|| \in L^1(\mu \otimes \mathbb{P})$ and  $\lambda(x,\omega)$ is almost surely constant and equal to 
\begin{align*}
\lambda := -\int_{M\times \Omega}\log||df_{\omega}(x)^{-1}||\,\,d\mu(x)d\mathbb{P}(\omega).
\end{align*}

\renewcommand\thesuposicao{(H3)} 
\begin{suposicao} \label{(C)} 
The Markov process associated to the random iterations $f^n_{\omega}$ admits an  ergodic stationary measure $\mu$ satisfying
\begin{align*}
\lambda = -\int_{M\times \Omega}\log||df_{\omega}(x)^{-1}||d\mu(x)d\mathbb{P}(\omega)> 0.
\end{align*}
\end{suposicao}
This condition implies that  for $\mathbb{P} \otimes \mu$ almost every $(x,\omega)$ all Lyapunov exponents  are positive, i.e. $f^n_{\omega}$ is non-uniformly expanding.

Our last assumption is used to obtain  large deviation estimates for the Lyapunov exponent $\lambda$, and ensures that orbits do not spend too much time too close to  the critical set.
\renewcommand\thesuposicao{(H4)} 
\begin{suposicao} \label{(D)} 
Let $m^*$ be the Lebesgue measure restricted to the support of the ergodic stationary measure $\mu$:
\begin{equation}\label{restricted_meas}
m^* = m_{|\mathrm{supp}(\mu)}. 
\end{equation}
Then, there exists a constant $C_1 >0$ such that,  for all $\delta$ small enough, the following  holds:
\begin{equation}\label{integral0}
 m^* \otimes \mathbb{P} \left\{(x,\omega) \colon \text{dist}(x,\mathcal{C}_{\omega_0}) \le \delta \right\} \le C_1 \delta.
\end{equation}
\end{suposicao} 
We will use  the two immediate consequences of \eqref{integral0}: there exist $C_2,C_3,a_0>0$ such that 
\begin{equation}\label{integral1}
\int_{\Omega} \int_{B_{\delta}(\mathcal{C}_{\omega_0})} -\log(\text{dist}(x,\mathcal{C}_{\omega_0})) d m^*(x) d \mathbb{P}(\omega)
\le  -C_2 \delta \log(\delta),
\end{equation}
and 
\begin{equation}\label{integral2}
\int_{\Omega}\int_{B_{\delta}(\mathcal{C_{\omega_0}}) } e^{-a \log(\text{dist}(x,\mathcal{C}_{\omega_0}))} d m^*(x) d \mathbb{P}(\omega_0)
\le   C_3 \delta^{1-a}, \qquad \forall a \le a_0.
\end{equation}
We note that \eqref{integral0} is satisfied whenever the critical set $\mathcal{C}_{\omega_0}$ admits a piecewise $C^{1+\alpha}$ parametrization with uniform bounds on $\omega_0$. 
\subsection{Random horseshoes}
We proceed to define our notion of a random horseshoe, following \cite{lamb2023horsehoes}.
\begin{definicao}\label{Horseshoe}
Let $B_0,B_1$ be two balls such that
$B_0 \cap B_1 = \emptyset $. Then, the pair $(B_0,B_1)$ is a  random horseshoe if the following holds: there exists 
a sequence of integer-valued random variables $\{n_k(\omega) \}_{k \ge 0}$ satisfying  almost surely 
\begin{equation}\label{e0}
\limsup_{k \to \infty} \frac{n_k(\omega)}{k} < C,
\end{equation}
with $C>0$ independent of $\omega$, 
\begin{equation}\label{e1} f^{n_{k+1}(\omega)-n_k(\omega)}_{\theta^{n_k(\omega)}(\omega)}(B_i) \supseteq B_0 \cup B_1 \qquad \forall i \in \{0,1 \},
\end{equation}
and there exist balls $J(k,\omega)_{i,j} \subset B_i$ such that 
\begin{equation}\label{idk}
f^{n_{k+1}(\omega)-n_k(\omega)}_{\theta^{n_k(\omega)}(\omega)}(J(k,\omega)_{i,j}) = B_j
\end{equation}
and
\begin{equation}\label{e2}
\sup_{x \in J(k,\omega)_{i,j}}\left|\left|\left(df^{n_{k+1}(\omega)-n_k(\omega )}_{\theta^{n_k(\omega)}(\omega)}(x)\right)^{-1}\right| \right| <\kappa^{-1} \qquad \forall  k \ge 0,
\end{equation}
for some $\kappa>1$ independent of $\omega$.
\end{definicao}
With this definition, if $(B_0,B_1)$ is a horseshoe and we consider the set $P(\omega)$, which consists of all points  $x$ such that  $f^{n_k(\omega)}_{\omega}(x) \in B_0 \cup B_1$, for all $k \ge 0$, i.e.
\begin{equation}\label{Dani}
P(\omega):= \{x \in M \colon f^{n_k(\omega)}_{\omega}(x) \in B_0 \cup B_1\, \forall k \ge 0 \},
\end{equation}
then, the set  $P(\omega)$ is almost surely hyperbolic and in one-to-one correspondence with the sequences of zeros and ones.

Our first main result is 
\begin{teorema}\label{thm1}
Let $f_{\omega}$ be the random map satisfying Hypotheses \ref{(A)}-\ref{(D)}. Then, the random horseshoes are dense, which means that every open subset of the support of the ergodic measure contains  balls $B_0,B_1$ such that $(B_0,B_1)$ is a  random horseshoe.
\end{teorema}
In fact, our construction will give that the sequence of random times $\{n_{k}(\omega)\}_{k \ge 0}$ associated to a horseshoe $(B_0,B_1)$ has independent increments and satisfies 
\begin{align*}
\lim_{k \to \infty} \frac{n_k(\omega)}{k} = \mathbb{E}[n_0]
\end{align*}
almost surely. In particular, the asymptotic frequency of return times is independent of $\omega$.
\subsection{Random Young towers}\label{gibbsmarov}
In this subsection we provide formal definition of a random Young tower (see \cite{ASENS_2002_4_35_1_77_0}) and then state our second main result establishing the existence of random Young towers with annealed exponential tails for systems satisfying Hypotheses \ref{(A)}-\ref{(D)}.

Essentially, a Young tower may be thought of as a discrete-time suspension of a return  map with the following properties.

\begin{definicao}\label{def:random_markov_map}
Given a reference set  \( \Delta \subset M \)
we say that a family of  piecewise smooth maps
\begin{align*}
F_{\omega} \colon  \Delta \to \Delta ,
\end{align*}
is a \emph{Random Gibbs-Markov Map} 
if 

\begin{enumerate}
    
     \item[(Q1)]\textbf{Markov Property:}  For \( \mathbb{P}\)-almost every \( \omega \in \Omega \), there exists a countable (mod 0)-partition $P(\omega)$ of $\Delta$, i.e $P(\omega)=\{P_j(\omega)\}_{j=1}^{\infty}$, such that,
      for each $j$, the map \( F_{\omega} \colon P_j(\omega) \to \Delta \) is a diffeomorphism,

    \item[(Q2)] \textbf{Expansion:} there exists a constant \( \kappa > 1 \), independent of \( \omega \), such that for each $j \in \mathbb{N}$
    \[
    \sup_{x \in P_j(\omega)} ||dF_\omega(x)^{-1}|| < \kappa^{-1},
    \]

    \item[(Q3)] \textbf{Bounded Distortion:} there exists a constant \( \bar C > 0 \), independent of \( \omega \), such that for each $j$ and \( x, y \in P_j(\omega) \),
    \[
    \left| \log \frac{|\det dF_\omega(x)|}{|\det dF_\omega(y)|} \right| \leq \bar C |F_\omega(x) - F_\omega(y)|^{\beta}.
    \]
\end{enumerate}
\end{definicao}
In the deterministic case, this definition goes back to Alves \cite{Alves,ALVES_VILARINHO_2013}, whilst its random counterpart can be found in \cite{Zhiku2015}.
\begin{definicao}\label{randomyoung}
We say that the random dynamical system $f_{\omega}$ carries a random Young tower if there is a reference set $\Delta$ with $\mu(\Delta)>0$  and a Gibbs-Markov map $F_{\omega}\colon \Delta \to \Delta $ such that 
\begin{align*}
F_{\omega}=  f^{R(x,\omega)}_{\omega}(x),
\end{align*}
where 
\begin{align*}
R \colon \Delta  \times  \Omega \to \mathbb{N}
\end{align*}
is  a measurable function (the return time function) which, at $\mathbb{P}$ almost every $\omega$ is constant on the elements of the partition $\mathcal{P}(\omega)$ associated to the Gibbs-Markov map $F_{\omega}$.
\end{definicao}

Our main result is the following.
\begin{teorema}\label{thm2}
Let $f_{\omega}$ be the random map satisfying Hypotheses \ref{(A)}-\ref{(D)}. Then, $f_{\omega}$ carries a  random Young tower with annealed exponential tails, in the sense that there exists $D_0 \ge 1$ and $\gamma>0 $ such that the return time function $R$ associated to the tower satisfies
\begin{equation}\label{annealedexponentialtail}
 m \otimes \mathbb{P} \left\{(x,\omega) \in \Delta \times \Omega \colon R(x,\omega)>n \right\} < D_0e^{-\gamma n}.
\end{equation}
\end{teorema}
Recall that, by ergodicity of $\mu$, $\mu(\Delta)>0$ implies that Lebesgue almost every orbit in the support of $\mu$ enters $\Delta$ and that the union of images of $\Delta$ is the support of $\mu$. 

Formally defined, the Young tower is a random dynamical system. If it is mixing, then \eqref{annealedexponentialtail} implies exponential decay of correlations and an almost sure invariance principle \cite{su2023random}, including central limit theorem. Several papers provide sufficient conditions for mixing of a Young tower for either random or sequential dynamical systems by means of a pathwise aperiodicity condition for the return times associated to the Young tower \cite{bahsoun2019quenched,su2023random,AlvesBahsounRuziboevVarandas2023,KorepanovLeppanen2021}. These conditions are not inherently satisfied in our setting, but we nevertheless expect the tower to be mixing. However, a detailed analysis of this lies beyond the scope of the present paper. 

\section{Large Deviation estimates}
\label{section3}
Given $r \in (0,1)$ and $\omega \in \Omega$, we define the $r$-truncated distance from $x$ to the critical set $\mathcal{C}_{\omega_0}$  as \begin{equation}\label{criticaldistance}
\mathrm{dist}_r(x,\mathcal{C}_{\omega_0}):= \begin{cases} 1
\qquad &\text{if}\,\, \mathrm{dist}(x,\mathcal{C}_{\omega_0})> r, \\
\mathrm{dist}(x,\mathcal{C}_{\omega_0}) \qquad &\text{otherwise}.
\end{cases}
\end{equation}
Let $m^*$ be the restriction of the Lebesgue measure $m$ to the support of the ergodic measure $\mu$ as introduced in \eqref{restricted_meas}.
Consider the following sums
\begin{equation}\label{birkhoff1}
S_n(x,\omega) :=  \sum_{i = 0}^n \log\left|\left|df_{\theta^i(\omega)}(f^{i}_{\omega}(x))^{-1}\right|\right|,
\end{equation}
and
\begin{equation}\label{birkhoff2}
Z_n(x,\omega,r) := - \sum_{i=0}^{n} \log \mathrm{dist}_r(f^i_{\omega}(x),\mathcal{C}_{\omega_i}).
\end{equation}
Observe that, by \eqref{BVestimate} and \eqref{integral1}
\begin{align*}
\int_{\Omega} \int_{M}-\log \mathrm{dist}_r(x,\mathcal{C}_{\omega_0}) d\mu(x) d\mathbb{P}(\omega)   \le  C_0 \int_{\Omega} \int_{M}-\log \mathrm{dist}_r(x,\mathcal{C}_{\omega_0})d m^*(x) d\mathbb{P}(\omega)<\infty,
\end{align*}
i.e  $-\log \mathrm{dist}_r(x,\mathcal{C}_{\omega_0}) \in L^1(\left(\mu \otimes \mathbb{P}\right))$. Similarly,  applying first  \eqref{c1} in Hypothesis \ref{(B)} and then again \eqref{integral1},\eqref{BVestimate}, we have $\log||df_{\omega}(x)^{-1} || \in L^1(\left(\mu\otimes\mathbb{P}\right))$, As a result, by Birkhoff ergodic theorem
\begin{align*}
\lim_{n \to \infty} \frac{Z_n(x,\omega,r)}{n} = \int_{\Omega} \int_{M}- \log \mathrm{dist}_r(x,\mathcal{C}_{\omega_0})d\mu(x) d\mathbb{P}(\omega) =O(-r\log(r)) \qquad  \mu \otimes \mathbb{P}  \,\,a.s.,
\end{align*}
 and, by Hypothesis \ref{(C)},
\begin{align*}
\lim_{n \to \infty} \frac{S_n(x,\omega)}{n} = \int_{\Omega} \int_{M}  \log||df_{\omega}(x)^{-1} || \,d\mu(x) d\mathbb{P}(\omega) = - \lambda <0 \qquad \mu \otimes \mathbb{P}  \,\,a.s.
\end{align*}
The goal of this section is to prove the following annealed  large deviation estimates for the above limits. 
\begin{proposicao}\label{devmultid}
Assume Hypotheses \ref{(A)}-\ref{(D)} hold. Then, there exist $C_4>0$, $\gamma >0$ such that, for all $r>0$ small enough and all $n \ge 0$
\begin{equation}\label{final_large}
m^* \otimes \mathbb{P}  \left\{S_n(x,\omega) > -\frac{\lambda}{2}n\;\; or  \;\; Z_n(x,\omega,r) > -\sqrt{r}\log(r)n \right\} < C_4 e^{-\gamma n}.
\end{equation}
\end{proposicao}

Define
\begin{equation}\label{tildeesse}
\tilde{S}_n(x,\omega) :=  \sum_{i = 1}^n \log\left|\left|df_{\theta^i(\omega)}(f^{i}_{\omega}(x))^{-1}\right|\right| 
\end{equation}
and 
\begin{equation}\label{tilde_zed}
\tilde{Z}_n(x,\omega,r) := \sum_{i=1}^{n} -\log \mathrm{dist}_r\left(f^i_{\omega}(x),\mathcal{C}_{\omega_i}\right).
\end{equation}
\begin{lema}\label{simp}
There exists $C_5\ge 0 $ and $z>0$ such that, for all $r>0$ small enough and all $n \ge 0$  
\begin{equation}\label{simpleeq}
 \mathbb{P}\left\{\tilde{S}_n(x,\omega) > -\frac{\lambda n}{4}\right\}
+ \mathbb{P}\left\{\tilde{Z}_n(x,\omega,r) > -\frac{\sqrt{r}\log(r)}{4}n \right\}< C_5e^{-z n} \qquad  \forall x \in \supp \mu
\end{equation}
\end{lema}

By  \eqref{c1} in Hypothesis \ref{(B)} and \eqref{integral0}  in Hypothesis \ref{(D)}, Lemma \ref{simp} implies Proposition \ref{devmultid} immediately.

The proof of this Lemma is deferred to the next subsections, and is an adaptation of the classical approach developed in  \cite{bougerol1985products}, where analytic perturbation theory is combined with the Nagaev method \cite{Nagaev} and the analyisis of the tilted operator, as defined in \eqref{tilted}, to prove large deviation estimates for bounded observables of  \(C^\alpha\)-mixing Markov chains.

In this work, we modify this scheme to treat the augmented chain \((f^n_{\omega}(x),\omega_n)\) and unbounded observables satisfying the exponential moment bounds stated in \eqref{integral2} (see Lemma \ref{LdevLemma1}). The main issue is that, for this Markov chain, the associated tilted operator, the main tool of the theory, is not even bounded  in $C^0$. We solve this problem by considering the spectrum of the tilted operator in the Banach space $\mathcal B(\tilde M)$, as defined in \eqref{space}, and not in the classical space \(C^0\). It is worth noting that large deviation estimates for unbounded observables were obtained, both for $C^0$-mixing Markov chains and for uniformly expanding maps, under weaker moment assumptions; see, for instance, \cite{Nagaev,herve2010nagaev}. However, this approach requires one to establish separately that the tilted operator has a spectral gap. In our setting of non-uniform expansion with unbounded observables, we obtain the spectral gap of the tilted operator through analytic perturbation theory, which requires the exponential moments assumption.

For the rest of Section \ref{section3}, in order to avoid unnecessary technicalities, we assume that
\[
\operatorname{supp}(\mu)=M
\qquad\text{and}\qquad
m^*=m.
\]
This assumption simplifies the computation drastically as it implies, by \cite[Theorem B]{CastroLambOliconMendezRasmussen2024}, that $\mu$ is mixing. 
If this assumption does not hold, the same arguments apply after restricting the dynamics to the relevant cycling components of $\supp \mu$ and applying the result to a suitable iteration of the random system $f_{\omega}$, which is mixing on such a component.


\subsection{The augmented Markov Chain}
Define
\begin{align*}
\tilde M:= M\times[-\sigma,\sigma]^k,
\end{align*}
and consider the Banach space
\begin{equation}\label{space}
\mathcal{B}(\tilde M):=
\left\{
\phi(x,\omega_0)\colon \tilde M\to\mathbb R
\,\,\vert\,
\sup_{x\in M}
\left(
\int |\phi(x,\omega_0)|^2d\mathbb P(\omega_0)
\right)^{1/2}<\infty
\right\},
\end{equation}
endowed with the norm
\begin{equation}\label{norm}
||\phi||_{\mathcal{B}}:= \sup_{x \in M } \left(\int |\phi(x,\omega_0)|^2d\mathbb{P}(\omega_0)\right)^{1/2}.
\end{equation}
It is easy to see that $C^0(\tilde{M}) \subset \mathcal{B}(\tilde M)$. We now consider the transfer operator 
\begin{align*}
&\mathcal{P} \colon \mathcal B(\tilde M) \to \mathcal{B}(\tilde M)\\
&\mathcal{P}(\phi)(x,\omega_0) =   \int \phi(g(x,\omega_0),\omega_1)d\mathbb{P}(\omega_1).
\end{align*}
Note that its second iterate is given by
\begin{equation}\label{kernelop}
\mathcal{P}^2(\phi)(x,\omega_0) = \int_{ g(u,[-\sigma,\sigma]^k)\times[-\sigma,\sigma]^k} \phi(z,\omega_1)  \left|\det(\frac{\partial g}{\partial \omega_1}(u,\bar \omega_1(z,u)))\right|^{-1}dm(z)d\mathbb{P}(\omega_1),
\end{equation}
where we denote $u = g(x,\omega_0)$ and $\bar \omega_1(z,u)$ is the solution to the equation $
z = g(u,\omega_1).$

Since, by Hypothesis \ref{(A)}, the map 
\begin{align*}
u \mapsto \left|\det(\frac{\partial g}{\partial \omega_1}(u,\bar \omega_1(z,u)))\right|^{-1}
\end{align*}
is continuous and bounded, uniformly for all $z$, we have
\begin{equation}\label{strongf}
\mathcal{P}^2(L^{\infty}(\tilde M)) \subset C^0(\tilde M),
\end{equation}
i.e $\mathcal{P}^2$ is strong-Feller.

\begin{proposicao}\label{pixi}
The operator $\mathcal{P}$ is eventually compact. It has a simple eigenvalue $\lambda = 1$ and the rest of the spectrum lies strictly inside the unit circle. 
\end{proposicao}
\begin{proof}
The strong-Feller property \eqref{strongf} implies, by \cite[Theorem 1.5.9]{Revuz1984}, that $\mathcal{P}^4$, restricted to $C^0(\tilde{M})$, is compact.  By \cite[Theorem B]{CastroLambOliconMendezRasmussen2024}  and the fact that $\tilde M$ is connected, the peripheral spectrum in $C^0(\tilde{M})$ consists only of the simple eigenvalue $1$.

By \eqref{kernelop}
\begin{align*}
|\mathcal{P}^2(\phi)(x,\omega_0)| &\le \left(\sup_{x \in M} \int |\phi(x,\omega_1)|d\mathbb{P}(\omega_1)\right) \,\left(\int_{ g(u,[-\sigma,\sigma]^k)} \left|\det(\frac{\partial g}{\partial \omega_1}(u,\bar \omega_1(z,u)))\right|^{-1}dm(z)\right)\\
&\le ||\phi||_{\mathcal B} \int_{[-\sigma,\sigma]^k} d\mathbb{P}(\bar \omega_1) = ||\phi||_{\mathcal B}.
\end{align*}
 Hence, $||\mathcal{P}^2(\phi)||_{\infty} \le ||\phi||_{\mathcal B}$ and, by \eqref{strongf},
 $
 \mathcal{P}^4(\mathcal{B}(\tilde{M})) \subset C^0(\tilde M).
 $ Compactnes of $\mathcal{P}^4$  restricted to $C^0(\tilde M)$ implies 
\[
\mathcal P^8:\mathcal B(\tilde M)\to \mathcal B(\tilde M)
\]
is compact, as required,  and that all the eigenfunctions are in $C^0(\tilde M)$. 
\end{proof}

\subsection{Analytic perturbation theory}\label{sub32}
Let $A(x,\omega_0)$ be an observable satisfying the following moment estimate: there exists $\tilde{\xi}>0$ such that for all $0\le  \xi\le \tilde \xi$ 
\begin{equation}\label{momentestimamte}
\int_{M\times [-\sigma,\sigma]^k} e^{\xi |A(x,\omega_0)|}dm(x)d\mathbb{P}(\omega_0)<\infty. 
\end{equation}

We define  the tilted operator $\mathcal P_{\xi,A}\colon \mathcal{B}(\tilde{M}) \to \mathcal{B}({\tilde{M}})$ associated to such $A$ as
\begin{equation}\label{tilted}  
\mathcal P_{\xi,A}(\phi)(x,\omega_0) =\int e^{\xi A\left(g(x,\omega_0),\omega_1\right)}\phi\left(g(x,\omega_0),\omega_1\right)d\mathbb{P}(\omega_1)
\end{equation}

In the following Lemma, we are going to characterize the spectrum of  $\mathcal{P}_{\xi,A}$, for $\xi$ small enough. 
\begin{lema}\label{LdevLemma1}
Let $A$ be an observable satisfying \eqref{momentestimamte}. Then, there exists $\xi_1>0$ such that for all $|\xi|<\xi_1$ the  operator $\mathcal{P}_{\xi,A}$ admits the decomposition  
\begin{equation}\label{specdec}
\mathcal{P}_{\xi,A}= \lambda_{\xi,A}\Pi_{\xi,A} + Q_{\xi,A},
\end{equation}
where $\Pi_{\xi,A}$ is the projection to the eigenfunction corresponding to the  simple eigenvalue $\lambda_{\xi,A}$ and the remainder operator  $Q_{\xi,A} = \mathcal{P}_{\xi,A}- \lambda_{\xi,A}\Pi_{\xi,A}$ satisfies
\begin{align*}
||Q_{\xi,A}^n|| \le Cr^n,
\end{align*}
for some $r<1$ independent of $\xi$.
Furthermore, the leading eigenvalue $\lambda_{\xi,A}$ satisfies 
\begin{equation}\label{lambdaxi}
\lambda_{\xi,A} = 1  + \xi \int A(x,\omega_0)d\mu(x) d\mathbb{P}(\omega_0)  +O(\xi^2).
\end{equation}
\end{lema}
\begin{proof}
It is sufficient to prove that $\mathcal{P}_{\xi,A}$ admits the analytic expansion in $\mathcal{B}(\tilde M)$
\begin{equation}\label{analytic_identuty}
\mathcal{P}_{\xi,A}(\phi) = \sum_{n=0}^{\infty} \frac{\xi^n}{n!}\mathcal{P}_{n,A}(\phi),
\end{equation}
where
\begin{equation}\label{quenne}
\mathcal{P}_{n,A}(\phi)(x,\omega_0) = \int \left(A(g(x,\omega_0),\omega_1)\right)^n \phi(g(x,\omega_0),\omega_1)d\mathbb{P}(\omega_1).
\end{equation}
Indeed, if the above holds true, the analytic perturbation argument in  \cite[Lemma 3.2]{bougerol1985products}, along with the spectral decomposition for $\mathcal{P}_0$ established in Proposition \ref{pixi}, implies the spectral decomposition in \eqref{specdec} and the  analyticity of $\lambda_{\xi,A}$ and of its unique eigenfunction $\rho_{\xi,A}$ satisfying $\int \rho_{\xi,A} d\mu d\mathbb{P} = 1$. In particular, $\lambda_{\xi,A}, \rho_{\xi,A}$ satisfy  
\begin{align*}
\lambda_{\xi,A}  = \int \mathcal{P}_{\xi,A}(\rho_{\xi,A})(x,\omega_0) d\mu(x)d \mathbb{P}(\omega_0) .
\end{align*}
By \eqref{analytic_identuty},\eqref{quenne}, and the invariance of $\mu \otimes \mathbb{P}$, we have 
\begin{align*}
\lambda_{\xi,A} &= 1 + \xi \int A (x,\omega_0)\rho_{\xi,A}(x,\omega_0)d\mu(x)d \mathbb{P}(\omega_0)  + O(\xi^2).
\end{align*}
Using the analyticity of $\rho_{\xi,A}$ and the fact that $\rho_0 = 1$, we get 
\begin{align*}
 \lambda_{\xi,A} &= 1 +\xi  \int A(x,\omega_0)d\mu(x)d \mathbb{P}(\omega_0)  + O(\xi^2),
\end{align*}
which shows \eqref{lambdaxi}.

In order to show \eqref{analytic_identuty}, we will show that the partial sums 
\begin{equation}\label{partialsums}
R_{N,\xi,A}(\phi) := \sum_{n=0}^{N-1} \frac{\xi^n}{n!}\mathcal{P}_{n,A}(\phi)
\end{equation}
are bounded linear operators in $\mathcal{B}(\tilde M)$ and converge to $\mathcal{P}_{\xi,A}$ in operator norm.
Indeed, by Cauchy-Schwarz inequality, and the definition of the norm $||\cdot ||_{\mathcal{B}}$ in \eqref{norm}, we have 
\begin{align*} 
||\mathcal{P}_{n,A}||_{\mathcal{B}\to \mathcal{B} } \le \sup_{x \in M} \left(\int \left| A(g(x,\omega_0),\omega_1) \right|^{2n} d\mathbb{P}(\omega_0) d\mathbb{P}(\omega_1)\right)^{\frac{1}{2}}.
\end{align*}
This is  finite by \eqref{momentestimamte}, because, changing the variable $\omega_0 \mapsto z= g(x,\omega_0),$ we obtain
\begin{align*}
\int \left| A(g(x,\omega_0),\omega_1) \right|^{2n} d\mathbb{P}(\omega_0) d\mathbb{P}(\omega_1) \le \frac{(2n)!}{(2\sigma)^k \tilde{\xi}^{2n}c_2}  \int e^{\tilde \xi|A(z,\omega_1)|}dm(z) d\mathbb{P}(\omega_1) < \infty,
\end{align*}
where $c_2$ is as in \eqref{boundsforinvmeas}. Hence, for all $n \ge 0$, $\mathcal{P}_{n,A}$ is a bounded linear  operator in $\mathcal{B}(\tilde M)$,  and so are the partial sums $R_N$ in \eqref{partialsums}.

Now we show that $R_{N,\xi,A}$ converges to $\mathcal{P}_{\xi,A}$ as $N \to \infty$ in operator norm. Indeed
\begin{align*}
\left|\mathcal{P}_{\xi,A}(\phi)(x,\omega_0) -R_{N,\xi,A}(\phi)(x,\omega_0) \right|&\le  \int \left|\sum_{n \ge N } \frac{\xi^n}{n!} A(g(x,\omega_0),\omega_1)^n \phi((g(x,\omega_0),\omega_1))\right|d\mathbb{P}(\omega_1) \\
&\le \left(\int \left|\sum_{n \ge N } \frac{\xi^n}{n!} A(g(x,\omega_0),\omega_1)^n\right|^2 d\mathbb{P}(\omega_1)\right)^{1/2}||\phi||_{\mathcal B},
\end{align*}
so that 
\begin{align*}
\left|\left|\mathcal{P}_{\xi,A} -R_{N,\xi,A} \right|\right|_{\mathcal{B} \to \mathcal{B}} &\le \sup_{x \in M}   \left(\int \left|\sum_{n \ge N } \frac{\xi^n}{n!} A(g(x,\omega_0),\omega_1)^n\right|^2 d\mathbb{P}(\omega_0)d\mathbb{P}(\omega_1)\right)^{1/2}\\
&\le \frac{1}{(2\sigma)^{k/2}c_2^{1/2}} \left(\int \left|\sum_{n \ge N } \frac{\xi^n}{n!} A(z,\omega_1)^n\right|^2 d\mathbb{P}(\omega_1)dm(z)\right)^{1/2}
\end{align*}
where $c_2$ is as in \eqref{boundsforinvmeas}.
Note that, almost surely in $(z,\omega_1)$,
\begin{align*}
\sum_{n \ge N } \frac{\xi^n}{n!} A(z,\omega_1)^n   \to 0 \qquad N \to \infty,
\end{align*}
and we have, uniformly in $N$,
\begin{align*}
\left|\sum_{n \ge N } \frac{\xi^n}{n!} A(z,\omega_1)^n\right|^2 \le e^{2\xi |A(z,\omega_1) |},
\end{align*}
whose integral is finite for all  $
\xi \le \frac{\tilde \xi}{2}$ due to \eqref{momentestimamte}.
 By the  dominated convergence theorem 
\begin{align*}
\lim_{N \to \infty} \left|\left|\mathcal{P}_{\xi,A} -R_{N,\xi,A} \right|\right|_{\mathcal{B} \to \mathcal{B}} = 0.
\end{align*}
\end{proof}
\begin{corolario}\label{ldevcorollary}
Given an observable $A$ satisfying \eqref{momentestimamte},
define
\begin{align*}
A_n(x,\omega) =\sum_{i=1}^n A(f^i_{\omega}(x),\omega_i),
\end{align*}
and
\begin{align*}
\lambda_A = \int A(x,\omega_0) d\mu(x)d\mathbb{P}(\omega_0).
\end{align*}
Then, there exists $C_A, \xi_A>0$ such that, for all $\varepsilon >0$
\begin{equation}\label{corollardev}
\mathbb{P}\left\{ A_n(x,\omega) > (\lambda_A+\varepsilon)n \right\} < C_A e^{-\frac{\xi_A }{2} \varepsilon^2
 n }, \qquad \forall n \ge 0,\forall x \in M.
\end{equation}
\end{corolario}
\begin{proof}
Observe that, for any $\xi>0$,
\begin{align*}
\mathbb{E}[e^{\xi A_n}] = \int_{[-\sigma,\sigma]^k}\mathcal{P}_{\xi,A}^n(1)(x,\omega_0)d\mathbb{P}(\omega_0),
\end{align*}
where $\mathcal{P}_{\xi,A}$ denotes the tilted operator. As a result, by Chernoff bound, if  $\varepsilon>0$ and $x \in M$
\begin{align*}
\mathbb{P}\left\{ A_{n}(x,\omega) > (\lambda_A+\varepsilon)n \right \}  \le e^{-(\lambda_A+\varepsilon)\xi  n} \int_{[-\sigma,\sigma]^k}\mathcal{P}^n_{\xi,A}(1)(x,\omega_0)d\mathbb{P}(\omega_0)
\le e^{-(\lambda_A+\varepsilon)\xi n } ||\mathcal{P}^n_{\xi,A}(1) ||_{\mathcal{B}(\tilde{M})},
\end{align*}
which, by \eqref{lambdaxi}, proves the claim.
\end{proof}

\begin{TheoremproofD}
Apply Corollary \ref{ldevcorollary} to the observables  $A_S = \log||df_{\omega}(x)^{-1}||$ and $A_Z= -\log(\text{dist}_r(x,\mathcal{C}_{\omega_0}))$, which both satisfy the exponential moment estimates \eqref{momentestimamte}, due to \eqref{c1} to \eqref{integral2}
\end{TheoremproofD}

\section{Existence of Young Times}
\label{section5}
The aim of this section is to define the Young times and establish their main properties. Essentially, they are times at which  balls with exponentially small radius expand uniformly to cover a reference set $\Delta$ (see Proposition \ref{DaniWhy}).

This section is divided into  three subsections:
\begin{itemize}
    \item In subsection \ref{sub1} we  revisit the theory of hyperbolic times developed in \cite[Chapter 6]{Alves}, as it is crucial in the construction of Young times;
    \item Then, in subsection \ref{sub2}, we define the Young times and establish their existence and their main properties;
    \item Finally, in subsection \ref{sub3}, we estimate the first moment of Young times. This is crucial in our horseshoe construction.
\end{itemize}

\subsection{Hyperbolic times and their main properties}\label{sub1}
Hyperbolic times were introduced in
\cite{Alves2000SRBMF} for deterministic systems and extended in \cite{AST_2003__286__25_0}  to the random case and have been  exploited extensively to construct Young towers \cite{ASENS_2002_4_35_1_77_0,doi:10.1142/S0219493718500272,castro2023random}, to develop thermodynamic formalism \cite{Stadlbauer2021-wu} and to prove stochastic stability \cite{ALVES_VILARINHO_2013}. In  this section, we recall the definition and the main properties closely following the lines of  \cite[Section 6]{Alves}.  

Let $f^n_{\omega}$ be as in \eqref{rdsiterations}. We start fixing a real number $b>0$ satisfying
\begin{equation}\label{alphaminusonebeta}
b < \frac{\min\{1,\alpha^{-1} \beta\}}{2},
\end{equation}
where $\alpha,\beta$ are as in Hypothesis \ref{(B)}. 
The following definition is taken from \cite[Definition 4.1]{ALVES_VILARINHO_2013}.
\begin{definicao}\label{hyperbolic times}
 Given $\eta \in (0,1)$ and $r>0$, we say that $n$ is a $(\eta,r)$-hyperbolic time for $(x,\omega) \in M \times \Omega$ if for every $0 \le k < n$
\begin{align*}
\prod_{i=k}^{n-1}\left|df_{\theta^i\omega} (f_{\omega}^i(x))^{-1}\right| <\eta^{(n-k)},\ \text{and } \mathrm{dist}_{r}(f^{k}_\omega(x), \mathcal{C}_{\theta^{k}\omega}) > \eta^{b (n-k)}, 
\end{align*}
where the truncated distance $\mathrm{dist}_{r}\left(x,\mathcal{C}_{\omega_0}\right)$ is defined as  in \eqref{criticaldistance}.
\end{definicao}
The following Proposition shows that, whenever $n$ is a hyperbolic time for $(x,\omega)$, $f^n_{\omega}$  uniformly expands a sufficiently small neighborhood of $x$ to a neighborhood of $f^n_{\omega}(x)$ of a fixed  radius $\delta$. This is the essential property of hyperbolic times needed to prove  Proposition \ref{DaniWhy}. 
\begin{proposicao}\cite[Proposition 4.9]{ALVES_VILARINHO_2013}\label{defVn}
Given $\eta ,r>0$, there exist $\tilde{\delta}=\tilde{\delta}(\eta,r)$ and  $C_6= C_6(\eta,r)>0,$ such that for every $(x,\omega)\in M \times \Omega$ for which $n\in\mathbb N$ is a $(\eta,r)$-hyperbolic time and for every $0<\delta<\tilde\delta$ the following hold:
\begin{enumerate}\label{HipBalls}
    \item[(i)] There exists an open neighborhood of $x$
    \begin{equation}\label{viennedeltaomegax}
     V_n^\delta(x,\omega) \subset B_{\delta \eta^{\frac{n}{2}}}(x) 
    \end{equation}
    such that $f^n_\omega$ maps $V_n^{\delta}(x,\omega)$ diffeomorphically to $B_{\delta} (f_\omega^n(x)).$
    \item[(ii)] for every $z,y\in V_n^\delta(x,\omega)$
    \begin{enumerate}
        \item[$a.$] $|f_{\omega}^{n-k}(z) -f_{\omega}^{n-k}(y)| \leq \eta^{\frac{k}{2}} |f^n_\omega(z) - f^n_\omega(y)|,$ for every $1\leq k \leq n$;
        \item[$b.$] $\displaystyle \left|\log \frac{|\det \d f^n_\omega (z)|}{|\det \d f^n_\omega (y)|} \right|\leq C_6  |f^n_\omega(z) - f^n_\omega(y)|^{\beta}$,
    \end{enumerate}
\end{enumerate}
where $\beta$ is as in \eqref{c2}.
\end{proposicao}

The following proposition shows that typical orbits have finite-time density of hyperbolic times. 
\begin{proposicao}\cite[Proposition~4.4]{ALVES_VILARINHO_2013}\label{hyperbolictimes}
Given $l,r>0$ and a sufficiently small $\varepsilon>0$,
there exists  $\theta > 0$ such that  if  
\[
-\sum_{i = 0}^n \log\left|\left|df_{\theta^i(\omega)}(f^{i}_{\omega}(x))^{-1}\right|\right| >nl
\quad \text{and} \quad
 \sum_{i=0}^{n} \log \mathrm{dist}_r(f^i_{\omega}(x),\mathcal{C}_{\omega_i}) >- \varepsilon n,
\]
for some $n$ and $(x,\omega)$, 
then $(x,\omega)$ has $(e^{-l/4}, r)$-hyperbolic times $1 \leq n_1 < \cdots < n_t \leq n$ with $t \geq \theta n$.
\end{proposicao}
Let $H_n(x,\omega,\eta,r)$  denote the set of $(\eta,r)$-hyperbolic times which are less or equal to $n$, i.e. 
\begin{equation}\label{haccaenne}
H_n(x,\omega,\eta,r):= \{i \in \{1,\dots,n \} \colon i \,\,\text{is a $(\eta,r)$ hyperbolic time for $(x,\omega)$} \}.
\end{equation}
 An immediate corollary of Proposition \ref{hyperbolictimes} and the annealed  large deviation estimates in Proposition \ref{devmultid} is that, for proper $\eta_0,r_0$,
 the $m^* \otimes \mathbb{P}$ probability that $H_n(x,\omega,\eta_0,r_0)$ does not have uniformly  positive density is exponentially small.
 \begin{corolario}\label{existence of the hyperbolic times}
There exists $\theta_0 \in (0,1)$, $C_7 \ge 1$, $\gamma_1 >0$, $\eta_0 \in (0,1)$ and $r_0>0$ such that for all $n \ge 0$
\begin{equation}\label{carinality of hyperbolic times}
m^* \otimes \mathbb{P} \left\{ |H_n(x,\omega,\eta_0,r_0)| \le \theta_0 n  \right\} \le C_7e^{-\gamma_1 n } 
\end{equation}
\end{corolario}
Note that, because of \eqref{carinality of hyperbolic times} and the Borel-Cantelli Lemma, almost every $(x,\omega)$ has infinitely many $(\eta_0,r_0)$-hyperbolic times, with asymptotic density at least $\theta_0$.

\subsection{Existence and positive density of Young times}\label{sub2}
Our main results here are  Proposition \ref{DaniWhy}, concerning the dynamical properties of Young times and Proposition \ref{densityYT}, that establishes the existence and positive density of Young times along with tail estimates.

Let us fix $\eta_0,r_0$ such that the conclusion of Corollary \ref{existence of the hyperbolic times} holds. Let $\delta_0$ be the  value obtained applying Proposition \ref{defVn}, with $\eta = \eta_0$ and $r=r_0$ and let $A$ be an open set in the support of $\mu$. The ergodicity, along with the null-measure of the critical set and the continuity of $x \to f_{\omega}(x)$, gives us the existence of some $(\bar x,\bar \omega) \in M \times \Omega $, an integer $N \in \mathbb{N}$, and positive numbers  $\delta_1, \delta_2> 0$ such that the forward orbit $\{f^i_{\bar \omega}(\bar x)\}_{i=0}^N$ is  $\frac{\delta_0}{2}$-dense, satisfies
\begin{align*}
B_{\delta_2}(f^N_{\bar \omega}(\bar x)) \Subset f^{N-j}_{\theta^j(\bar \omega)}(B_{\delta_1}(f^j_{\bar \omega}(\bar x))) \Subset A\qquad \forall j \le N,
\end{align*}
 and stays  bounded away from the critical set. Define the reference set
\begin{equation}\label{referenceset}
\Delta:= B_{\delta_2}(f^N_{\bar \omega}(\bar x)).
\end{equation}
Then, since $\Delta \subset A$, the continuity of the map $\omega \to f_{\omega}(x)$, along with the upper semi-continuity of $\omega_0 \to \mathcal{C}_{\omega_0}$, as stated in Hypotheses \ref{(A)}-\ref{(B)}, imply that there exist  constants  $\iota,\rho>0$ such that 
\begin{equation}\label{lowerboundtransitivityannealed}
\mathbb{P}(\mathcal{E}_{\Delta}(B))>\rho \qquad \text{for every ball}\,\, B \subset \supp(\mu) \colon r(B) \ge \delta_0,
\end{equation}
where $\mathcal{E}_{\Delta}(B)$ is defined as the set of noise realizations for which the image of $B$  covers $\Delta$ in at most $N$ iterates while staying away from the critical set:
\begin{equation}\label{DaniTVB}
\mathcal{E}_{\Delta}(B):=  \left\{\exists i \le N \colon \exists \tilde{B} \subset B \colon f^i_{\omega}(\tilde{B})\,\, \Supset \,\, \Delta\, \text{and}\, \inf_{j \le i} \inf_{x \in \tilde{B}} d\left(f^j_{\omega}(x),\mathcal{C}_{\omega_j}\right)>\iota \right\}.
\end{equation}
From now on fix the values of $\delta_0,\rho,\iota$ and $N$ as above. Consider now the sequence  $\{\tau^N_i:M \times \Omega \to\mathbb N\cup \{\infty\}\}_{i\in\mathbb N},$ of \emph{$N$-sparse hyperbolic times}: $$ \tau^N_i(x,\omega) :=\begin{cases}
\min\{n; \ n\ \text{is a $(\eta_0,r_0)$-hyperbolic time for}(x,\omega)\}&, \ \text{if}\ i=1, \\
\min\{n; \ n >  \tau^N_{i-1}(x,\omega)+N\ \text{is a $(\eta_0,r_0)$-hyperbolic time for }(x,\omega)\}&,\ \text{if }i>1.
\end{cases} $$
Let  $\mathcal{T}_n^N(x,\omega)$ denote the set of $N$-sparse hyperbolic times which are less or equal to $n$, i.e. 
\begin{equation}\label{tienne}
\mathcal{T}^N_n(x,\omega):= \{1,\dots, n \} \cap \{\tau^N_i(x,\omega)\}_{i \ge 1}.
\end{equation}

It is clear that
\begin{equation}\label{posdensitynsparse}
|\mathcal{T}^N_n(x,\omega)| \ge \frac{|H_n(x,\omega,\eta_0,r_0)|}{N+1}-1.
\end{equation}

Now we are ready to define the concept of Young times.
\begin{definicao}\cite[Definition 3.6]{castro2023random}\label{youngtimes}
Given a reference set $\Delta$, $\delta_0,N$ as above and $n \in \mathbb{N}$ consider the set
\begin{equation}\label{youngtimeseq}
Y^{\Delta}_n(x,\omega) := \left\{i\in \mathcal{T}^N_n(x,\omega) \vert  \ \theta^{i}\omega \in \mathcal E_{\Delta}\left(B_{\delta_0}\left(f^i_\omega(x)\right)\right)\right\}.
\end{equation} 
We say that $i$ is an $(x,\omega,\Delta)$-Young time if there exists $n\in\mathbb N$ such that $i\in Y^{\Delta}_n(x,\omega)$.
\label{YoungTijmen}
\end{definicao} 
Note that, by construction, $(x,\omega)$-Young times are measurable in $\omega$ and upper-semicontinuous in $x$.

The importance of Young times is due to the fact that if $n$ is a Young time, then there exists $k \le N$ such that $f^{n+k}_{\omega}$  uniformly expands a small neighborhood of $x$ into covering $\Delta$  with uniform distortion bounds, as stated in the following proposition.
\begin{proposicao}\cite[Lemma 4.1]{castro2023random}\label{DaniWhy}
Suppose $n$ is an $(x,\omega,\Delta)$-Young time.
Then, there exists $k \le N$ and an open set, 
$$\tilde{B} \subset B_{\delta_0\eta_0^{\frac{n}{2}}}\left(x\right),$$
diffeomorphic to a ball,  such that $f^{n+k}_{\omega}(\tilde{B})= \Delta$, where $\Delta$ is defined in \eqref{referenceset}.
Furthermore, there exists a constant $ C_8 \ge 1$ such that
\begin{equation}\label{coolexp}
||\left(df^{n+k}_{\omega}(x)\right)^{-1}||< C_8\eta_0^{\frac{n+k}{2}} \qquad \forall x \in \tilde{B},
\end{equation}
and 
\begin{equation}\label{distortionboundsdelta}
\left|\log\left|\frac{\det df^{n+k}_{\omega}(x)}{\det df^{n+k}_{\omega}(y)} \right| \right|\le C_8 \qquad \forall x,y \in \tilde{B}.
\end{equation}
\end{proposicao}

The following Proposition gives tail estimates for the density of the Young times.
\begin{proposicao}\label{densityYT}
Given an   open set $A \subset M$, let $\Delta \subset A$ be the reference set given by equation \eqref{referenceset}. Then, there exists $C_9 \ge 1,$ $\gamma_2>0$ and $\theta_1 \in (0,1)$, such that for all $n \ge 1$
\begin{equation}\label{LargeDeviationsYoungTimes}
m^* \otimes \mathbb{P} \left\{|Y^{\Delta}_n(x,\omega)| \le \theta_1 n \right\} \le C_9e^{-\gamma_2 n}.
\end{equation}
\end{proposicao}
\begin{proof}
Let 
\begin{equation}\label{choiceofu}
\theta_2:= \frac{\theta_0}{2N+2},
\end{equation}
where $\theta_0$ is as in \eqref{carinality of hyperbolic times}.
By \eqref{posdensitynsparse} and choosing $n$ large enough, we have 
\begin{equation}\label{inclusiondenseN}
\left\{|\mathcal{T}^N_n(x,\omega)| \le \theta_2 n \right\} \subset \left\{|H_n(x,\omega,\eta_0,r_0)| \le \theta_0 n \right\}. 
\end{equation}
Therefore, by  Corollary \ref{existence of the hyperbolic times}, there exists $C_7 \ge 1$ and $\gamma_1>0$ such that, for all $n \ge 1$
\begin{equation}\label{large_dev_n_sparse}
m^* \otimes \mathbb{P} \left\{|\mathcal{T}^N_n(x,\omega)| \le \theta_2 n \right\} \le C_7e^{-\gamma_1 n }.
\end{equation}
Hence, for any $ t>0$ and $n \in \mathbb{N}$, one has 
\begin{align*}
m^* \otimes \mathbb{P}\left\{|Y^{\Delta}_n(x,\omega)| \leq t n \right\} &\le m^* \otimes \mathbb{P} \left\{|Y^{\Delta}_n(x,\omega)| \leq t  n, |\mathcal{T}^N_n(x,\omega)| \leq \theta_2 n\right\} \\ &+   m^* \otimes \mathbb P\left\{\ |Y^{\Delta}_n(x,\omega)| \leq t  n , |\mathcal{T}^N_n(x,\omega)| \geq \theta_2 n \right\} \\
&\le m^* \otimes \mathbb{P} \left\{|Y^{\Delta}_n(x,\omega)| \leq t  n, |\mathcal{T}^N_n(x,\omega)| \geq \theta_2 n\right\} +  C_7e^{-\gamma_1 n }.
\end{align*}
We have
\begin{align*}
m^* \otimes \mathbb{P}\left\{|Y^{\Delta}_n(x,\omega)| \leq t  n, |\mathcal{T}^N_n(x,\omega)|\geq \theta_2 n\right\} = \int_{\supp \mu} \mathbb{P}\left\{|Y^{\Delta}_n(x,\omega)| \leq t  n, |\mathcal{T}^N_n(x,\omega)| \geq \theta_2 n\right\} dm(x)\\
\end{align*}
which means that, in order to prove \eqref{LargeDeviationsYoungTimes}, it is sufficient to show that there exists $\bar D \ge 1$ and $\bar \gamma>0$ such that
\begin{equation}\label{sufficientyoungtimes}
\sup_{x \in \supp \mu} \mathbb{P}\left\{|Y^{\Delta}_n(x,\omega)| \leq t  n, |\mathcal{T}^N_n(x,\omega)| \geq \theta_2 n\right\}  \le \bar De^{-\bar \gamma n },
\end{equation}
for $t>0$ sufficiently small.
This is done exactly as in  \cite[Theorem 3.7]{castro2023random}.
\end{proof}
\subsection{Young times for balls}\label{sub3}
For any small ball $B$, let us define 
\begin{equation}\label{enneo}
n_0(B) := 2\frac{\log(7r(B)/8\delta_0)}{\log(\eta_0)},
\end{equation}
where $\delta_0,\eta_0$ are as in Proposition \ref{DaniWhy} and $r(B)$ denotes the radius of $B$.
\begin{definicao}\label{YoungTimesForballs}
 We say that $n$ is a $(\omega,B,\Delta)$-Young time if $n$ is an $(x,\omega,\Delta)$-Young time for a point $x$  at distance less or equal to $r(B)/8$ from the center of $B$ and $n > n_0(B)$.
\end{definicao}
Observe that, by Proposition \ref{DaniWhy} if $n$ is an $(\omega,B,\Delta)$-Young time, then there exists an open set $\tilde{B}$ fully contained in $B$ such that 
$f^{n+k}_{\omega}(\tilde{B})$ uniformly expands onto $\Delta$.

Given a ball $B$, define
\begin{equation}\label{firstYoungtimeforballs}
h(\omega,B,\Delta):= \min \{n \ge 1 \colon n \,\,\text{is an $(\omega,B,\Delta)$-Young time}\}.  
\end{equation}

In other words, $h(\omega,B,\Delta)$ is the minimum of the $(x,\omega,\Delta)$-Young times, that are larger than $n_0(B)$, across all points $x$ in the ball of radius $\frac{r(B)}{8}$ around the center of $B$. The measurability of $h$ follows from the upper-semicontinuous dependence of the Young times on $x$, so the minimum is taken on a countable dense subset of that ball.

In order to construct horseshoes, we need to estimate the tail $\mathbb{P}\left\{ h(\omega,B,\Delta )> k \right\}$.
\begin{proposicao}\label{expect}
Let $B$ be a ball contained in $\supp \mu$. If $\Delta$ satisfies Proposition \ref{densityYT}, then there exist constants $D_1,D_2 >0$ such that 
\begin{equation}\label{Firstmoment}
\mathbb{E}h(\omega,B,\Delta) \le D_1\log(r(B)^{-1}) +D_2.
\end{equation}
\end{proposicao}
\begin{proof}
For any $r>0$, let $V(r)$ 
denote the infimum of the volume  of all the balls of radius r in $M$.
Let  
\begin{equation}\label{kappa0}
k_0 := - \frac{2\log\left({V(r(B)/8)}\right)}{\gamma_2},
\end{equation}
with $\gamma_2$ as in \eqref{LargeDeviationsYoungTimes}
Consider, for $k \ge 1$, the  sequence of random sets
\begin{equation}\label{rarevent}
E_k(\omega):= \{x \in \supp\, \mu \colon |Y^{\Delta}_k(x,\omega)| < \theta_1 k \}.
\end{equation} 
Note that 
\begin{equation}\label{switch}
\mathbb{P}\left\{ m(E_k(\omega))>  V(r(B)/8) \right\} \le C_9e^{-\gamma_2\frac{k}{2} },\qquad \forall k \ge k_0,
\end{equation}
with $C_9$ as in \eqref{LargeDeviationsYoungTimes}.
Indeed, if there exists $k_1\ge k_0$ such that 
\begin{align*}
\mathbb{P}\left\{ m(E_{k_1}(\omega))>  V(r(B)/8)  \right\} > C_9e^{-\gamma_2\frac{k_1}{2}},
\end{align*}
then,  
\begin{align*}
\int_{\Omega} m(E_{k_1}(\omega)) d\mathbb{P}(\omega) > e^{-\gamma_2\frac{k_1}{2}} \mathbb{P}\left\{  m(E_{k_1}(\omega)) > V(r(B)/8) \right\}> C_9e^{-\gamma_2 k_1 },
\end{align*}
since $$V(r(B)/8)>e^{-\gamma_2 k/2}\qquad \forall k \ge k_0,$$ but this contradicts \eqref{LargeDeviationsYoungTimes}.

Let 
$$
K := \max\left\{n_0(B)/\theta_1,k_0 \right\},
$$
with $n_0$ as in \eqref{enneo}, $k_0$ as in \eqref{kappa0} and $\theta_1$ as in \eqref{LargeDeviationsYoungTimes}. 
Note that, because of the definition of $K$, there exists $\tilde{D}_1,\tilde{D}_2\ge 0$ such that
\begin{align*}
K \le \tilde{D}_1 \log(r(B)^{-1})+\tilde{D}_2.
\end{align*}
We claim that for any $k \ge K$
\begin{equation}\label{helpful_inclusion}
\left\{\omega \in \Omega \colon h(\omega,B,\Delta)>k \right\} \subset \{ \omega \in \Omega \colon  m(E_k(\omega)) > V\left(r(B)/8\right) \}.
\end{equation}
Indeed, suppose $m(E_k(\omega))< V\left(r(B)/8\right)$. Then there exists at least one point $x$ in the ball of radius $r(B)/8$ around the center of $B$ such that $x \notin E_k(\omega)$.
By the definition of $E_k$ in \eqref{rarevent}, there exists at least $\theta_1 k$ of $(x,\omega,\Delta)$-Young times less or equal than $k$. The largest of these Young times is at least  $\theta_1 k \ge \theta_1 K \ge n_0$. As a result, this $(x,\omega,\Delta)$-Young time is an $(\omega,B,\Delta)$-Young time by Definition \ref{YoungTimesForballs}, which proves \eqref{helpful_inclusion}.

Therefore, by \eqref{switch}-\eqref{helpful_inclusion}
\begin{equation}\label{averageestimate}
\mathbb{E}[h(\omega,B,\Delta)]
\le \sum_{k \le K} 1 + C_9\sum_{k \ge K} e^{-\gamma_2 k/2} \le  -D_1 \log(r(B)) + D_2,
\end{equation}
for some $D_1,D_2 >0$, which concludes the proof.
\end{proof}
\section{Existence of Horseshoes}
\label{section6}
The aim of this section
is  to prove  Theorem \ref{thm1} about the existence and abundance of random horseshoes. We start with a  few general statements that will be used later.
\subsection{Preliminaries}
The first Proposition ensures that the Birkhoff averages of IID random variables with infinite expectation pointwise explode almost surely. 
\begin{proposicao}\cite[Theorem 2.4.5]{Durrett2019}\label{infaverage}
Let $\{X_m \}_{m\ge 0}$ be a sequence of non-negative iid random variables  such that $E[X_0]= \infty$. Then
\begin{align*}
\liminf_{m \to \infty} \frac{X_0+\dots+X_{m-1}}{m}= \infty, \qquad a.s.
\end{align*}
\end{proposicao}

Now suppose we have an increasing sequence $\{n_j\}_{j\ge 0}$. The following two propositions establish a relationship between the asymptotic behavior of $\left\{\frac{n_j}{j}\right\}_{j \ge 0}$ and the density of $\{n_j\}_{j\ge 0}$.

\begin{proposicao}\label{equivalence1}
Let $\{n_j\}_{j \ge 0}$ be an increasing sequence. Then  
\begin{equation}\label{positived1}
\limsup_{k \to \infty} \frac{n_k}{k} <  \bar \alpha  \implies \liminf_{k \to \infty } \frac{|\{ j \colon n_j \le k \}|}{k} >  \frac{1}{\bar \alpha},
\end{equation}
and 
\begin{equation}\label{positived2}
\liminf_{k \to \infty } \frac{|\{ j \colon n_j \le k \}|}{k} > \bar  \beta \implies \limsup_{k \to \infty} \frac{n_k}{k}<  \frac{1}{\bar \beta}.
\end{equation}
\end{proposicao}
\begin{proof}
First we prove \eqref{positived1}. The assumption
$\limsup_{k \to \infty} \frac{n_k}{k} <  \bar \alpha$
implies that there exists $k_{0}>0$ such that for all $k >k_0$
\begin{align*}
n_k <  \bar{\alpha} k.
\end{align*}
Note that, since the sequence  is increasing, this gives  
\begin{align*}
\left| j \colon n_{j} \le   \bar \alpha k\right| \ge k,
\end{align*}
hence
\begin{align*}
\liminf_{k \to \infty} \frac{\left| j \colon n_{j} \le  k\right|}{k}  = \liminf_{k \to \infty} \frac{\left| j \colon n_{j} \le  \bar \alpha k\right|}{\bar \alpha k}  > \frac{1}{\bar \alpha}.
\end{align*}
The proof of \eqref{positived2} is analogous.
\end{proof}
\begin{proposicao}\label{equivalence2}
Let $\{n_j\}_{j \ge 0}$ be an increasing sequence. Suppose that 
\begin{equation}\label{firstitem}
\limsup_{k \to \infty} \frac{|\{j\colon n_j \le  k \}|}{k} >\tilde  \beta.
\end{equation}
Then 
\begin{equation}\label{seconditem}
\liminf_{k \to \infty} \frac{n_k}{k} < \frac{1}{\tilde \beta}.
\end{equation}
\end{proposicao}
\begin{proof}
Suppose \eqref{firstitem} holds. Then, there exists an increasing sequence $k_r \to \infty$  such that
\begin{align*}
\liminf_{r \to \infty} \frac{|\{j \colon n_j \le k_r |\}}{k_r}>\tilde \beta.
\end{align*}
Then, there exists $r_0>0$ such that for all $r> r_0$ 
\begin{align*}
|\left\{j \colon n_j \le k_r \right\}| > \tilde \beta k_r.
\end{align*}
Then, for $r > r_0$,  $n_{\lfloor\beta k_r\rfloor} \le k_r$, which  implies 
\begin{align*}
\limsup_{r \to \infty} \frac{n_{\lfloor\tilde \beta k_r\rfloor}}{\lfloor \tilde  \beta k_r \rfloor} < \frac{1}{\tilde \beta}.
\end{align*}
\end{proof}
\subsection{Horseshoe construction}
In this subsection, we establish density of random horseshoes and conclude the proof of Theorem \ref{thm1}. 


\begin{proof}[Proof of Theorem~\ref{thm1}]

Let $A$ be an arbitrary open set. Let us choose the reference set $\Delta \subset A$ as in \eqref{referenceset}. For a large enough  $L$, whose precise value will be determined later, let $B_1,\dots,B_L$ be a sequence of disjoint balls of radius between $\frac{r(\Delta)}{4L}$ and $\frac{r(\Delta)}{L}$.

Let $h(\omega,B_i,\Delta)$ be the first Young time for $(\omega,B_i, \Delta)$. Define the \emph{return time}
\begin{equation}\label{tildeacca}
\tilde{h}(\omega,B_i,\Delta) = h(\omega,B_i,\Delta)+ k,
\end{equation}
where $k \le N$ is the smallest number, given by  Proposition \ref{DaniWhy}, for which there exists an open set $\tilde{B}_i$ with $\tilde{B}_i \Subset B_i$ such that
\begin{align*}
f^{h+k}_{\omega}(\tilde{B}_i) = \Delta, 
\end{align*}
and equations \eqref{coolexp}-\eqref{distortionboundsdelta} hold with $n = h$. Note that $\tilde{h}$ is a stopping time, i.e. it depends only on the past. Also note that, by choosing $L$ large enough, we can assume that $r(B_i)$ is small enough, so that equation \eqref{coolexp} gives us the uniform expansion for  $f_{\omega}^{h+k}$ on $\tilde{B}_i$ (note that $h \ge n_0(B_i)$ for all $i$ with $n_0$ as in \eqref{enneo}, by Definition \ref{YoungTimesForballs} of the first Young time).

Thus, for the balls $B_1,\dots,B_L$ we can define inductively the sequence of \emph{expanding return times},  by the rule 
$$\tilde{h}_0(\omega,B_i)=0,\qquad\tilde{h}_1(\omega,B_i)=\tilde{h}(\omega,B_i,\Delta),
$$
and
\begin{align*}
 \tilde{h}_{k+1}(\omega,B_i) = \tilde{h}(\theta^{  \tilde{h}_k(\omega,B_i)}(\omega),B_i,\Delta)+\tilde{h}_{k}(\omega,B_i).   
\end{align*}

Note that, by construction, for all $k\ge0$ and for all $i$, there exists an open set  $\tilde{B}_{k,i}$ such that
\begin{align*}
f^{\tilde h_k(\omega,B_i)}_{\omega}(\tilde{B}_{k,i}) = \Delta \supset \cup_j B_j,
\end{align*}
and 
$f^{\tilde h_k(\omega,B_i)}_{\omega}$  is uniformly expanding on $\tilde{B}_{k,i}$.

Since
\begin{align*}
f^{\tilde{h}_{k+1}(\omega,B_i)-\tilde{h}_{k}(\omega,B_i)}_{\theta^{  \tilde{h}_k(\omega,B_i)}(\omega)}(B_i) \supset B_i,
\end{align*}
and since the noise is iid, $\{\tilde{h}_k(\omega,B_i) \}_{k \ge 0}$ has independent and identically distributed increments. As a result, by Proposition \ref{expect}, since $\tilde h- h$ is uniformly bounded (see \eqref{tildeacca}), there exists $D_1,D_2>0$ such that
\begin{equation}\label{application}
\lim_{k \to \infty} \frac{\tilde{h}_k(\omega,B_i)}{k} = \mathbb{E}\tilde{h}(\omega,B_i) <  D_1\log\left(\frac{4L}{r(\Delta)}\right) +D_2 \qquad \mathbb{P}\,\,a.s.
\end{equation}

Now for every pair of balls  $B_i,B_j$ we define the \emph{simultaneous expanding return times} inductively: \begin{align*}
h_1^{i,j}(\omega):= \min\{h \ge 0 \colon  \tilde{h}_p(\omega,B_i)= \tilde{h}_q(\omega,B_j) =h\,\,\,\text{for some $p,q \in \mathbb{N}$} \},
\end{align*}
and 
\begin{align*}
h^{i,j}_{k+1}(\omega):= h^{i,j}_{k}(\omega)+h_1^{i,j}(\theta^{h^{i,j}_{k}(\omega)}(\omega)),
\end{align*}
so that
\begin{equation}\label{llnevents}
h^{i,j}_{k}(\omega) =h_1^{i,j}(\omega)+ \sum_{l=1}^{k-1} h_1^{i,j}(\theta^{h^{i,j}_{l}(\omega)}(\omega)), \qquad \forall k \ge 1.
\end{equation}
By this definition, at $h=h^{i,j}_{1}(\omega)$ the map $f^{h}_{\omega}$   simultaneously expands $B_i$ and $B_j$ onto $B_i \cup B_j$. Therefore, in order to prove the theorem, it is enough to show that there exists $i,j$ such that the corresponding sequence of return times $n_k := h_k^{i,j}$ satisfies the positive density condition  \eqref{e0}:
\begin{equation}\label{horsh}
\limsup_{k \to \infty}\frac{h^{i,j}_k(\omega)}{k}<C,\qquad \mathbb{P} \,\,a.s.
\end{equation}

Let 
\begin{equation}\label{viemm}
V_L=\frac{1}{\left( D_1\log\left(\frac{4L}{r(\Delta)}\right)+D_2\right)},
\end{equation}
and  take 
\begin{equation}\label{twozedem}
C= C_L:= \frac{L^2-L}{2(LV_L-1)}.
\end{equation}
\begin{proposicao}\label{positive_deense}
Let $L$ be large enough so that $LV_L>1$. Then, there exist $i,j$, independent of $\omega$, such that 
\begin{equation}\label{notlimsupposprob}
\liminf _{k \to \infty} \frac{h^{i,j}_k(\omega)}{k}< C_L,
\end{equation}
with positive probability.
\end{proposicao}
\begin{proof} 
Let us first show that, for almost every $\omega \in \Omega$, there exist $i(\omega),j(\omega)$ such that
\begin{equation}\label{notlimsup}
\liminf _{k \to \infty} \frac{h^{i(\omega),j(\omega)}_k(\omega)}{k}< C_L.
\end{equation}
By Proposition \ref{equivalence2}, it is sufficient to prove that, for almost every $\omega \in \Omega$, there exists $i(\omega),j(\omega)$ such that the intersection of corresponding sequences of expanding return times $\{ \tilde{h}_k(\omega,B_{i(\omega)})\}_{k \ge 0}$ and $\{ \tilde{h}_k(\omega,B_{j(\omega)})\}_{k \ge 0}$ has upper density larger than $C_L^{-1}$:
\begin{equation}\label{veryimportant}
\limsup_{m \to \infty}\frac{\left|\tilde H^{m}_{j(\omega)}(\omega) \cap \tilde H^{m}_{i(\omega)}(\omega) \right|}{m} > \frac{1}{C_L},
\end{equation}
where we denoted 
\begin{equation}\label{family}
\tilde H_i^{m}(\omega) := \{ \tilde{h}_k(\omega,B_{i})\}_{k \ge 0} \cap \{ 1,\dots, m\}.
\end{equation}
By \eqref{application} and Proposition \ref{equivalence1} we see that
\begin{equation}\label{liminf}
\liminf_{m \to \infty} \frac{|\tilde H^m_i(\omega)|}{m} > V_L.
\end{equation}
 By the inclusion-exclusion principle, for every $m$, we have
\begin{equation}\label{inclusion-exclusion}
 m \ge \left| \bigcup_{j=1}^L \tilde H^{m}_j(\omega) \right| \ge \sum_{j = 1}^L | \tilde H^{m}_j(\omega) | - \sum_{1 \le j < i \le L} \left|\tilde H^{m}_j(\omega) \cap \tilde H^{m}_i(\omega) \right|.
\end{equation}
Hence, from \eqref{liminf}, we see that  for almost every $\omega \in \Omega$ there exists $m_0(\omega)$ such that for all $m \ge m_0(\omega)$
\begin{align*}
\sum_{1 \le j < i \le L} \left|\tilde H^{m}_j(\omega) \cap \tilde H^{m}_i(\omega) \right| &\ge \sum_{j = 1}^L | \tilde H^{m}_j(\omega) | - m \\ 
&\ge (LV_L-1)m.
\end{align*}
Since there are $\frac{L^2-L}{2}$ pairs of $(i,j)$ in the sum above, for every $m > m_0(\omega)$ there exist indices $\tilde{i}(m,\omega)$ and $\tilde{j}(m,\omega)$ such that
\begin{equation}\label{interections}
\left|\tilde H^{m}_{\tilde{j}(m,\omega)}(\omega)\cap \tilde H^{m}_{\tilde{i}(m,\omega)}(\omega) \right|> \frac{m}{C_L}.
\end{equation}
Since the number of possible pairs $(i,j)$ is finite, there exist 
$i(\omega)$ and $j(\omega)$ such that
\begin{equation}\label{intersection}
\left|\tilde H^{m}_{j(\omega)}(\omega) \cap \tilde H^{m}_{i(\omega)}(\omega) \right| >m \frac{1}{C_L} \qquad \text{for infinitely many $m$},
\end{equation}
which shows \eqref{veryimportant}, hence \eqref{notlimsup}.

This means that
\begin{align*}
\mathbb{P}\left(\cup_{1\le i<j\le M}\left\{\omega \in \Omega \colon \liminf_{k\to\infty}\frac{h_k^{i,j}(\omega)}{k}<C_L \right\}\right)=1.
\end{align*}
Hence, there exists $i,j$ such that
\begin{align*}\label{liminfcool}
\mathbb{P}\left(\liminf_{k \to \infty} \frac{h_k^{i,j}(\omega) }{k} < C_L\right)>0.
\end{align*}
\end{proof}
In order to conclude the proof of Theorem \ref{thm1}, we recall that, the sequence $\{h^{i,j}_k(\omega) \}_{k \ge 0}$ has independent and identically distributed increments. By Proposition \ref{infaverage}, formula \eqref{notlimsupposprob}  implies that   $\mathbb{E}[h^{i,j}_1]$ is finite. Therefore, applying the Law of Large Numbers to \eqref{llnevents}, we have
\begin{align*}
\lim_{k \to \infty} \frac{h_k^{i,j}(\omega) }{k} = \liminf _{k \to \infty} \frac{h_k^{i,j}(\omega) }{k}, \qquad \mathbb{P}\, a.s.
\end{align*}
which proves \eqref{horsh}, by Proposition \ref{positive_deense}.
\end{proof}

\section{Young tower with annealed exponential tails}
\label{section7}
In this section we prove Theorem \ref{thm2}. The construction is based on the positive density of the Young times and the techniques presented in  \cite[Chapter 5]{Alves}, adapted in \cite{castro2023random} to the random setting for predominantly expanding random circle endomorphisms.

\begin{proof}[Proof of Theorem~\ref{thm2}]
Let us fix an open set $A$ and the set $\Delta \subset A$ as in the statement of Proposition \ref{densityYT}. Because of the tail estimate for Young times in \eqref{LargeDeviationsYoungTimes},  Lemma $4.1$ in \cite{castro2023random} ensures that, for $\mathbb{P}$-almost every noise realization, the hypotheses of the inducing construction in \cite[Chapter 5.3.2]{Alves} (see \cite[Corollary 5.9]{Alves}) are satisfied. We  therefore apply this construction pathwise, -as essentially done in Step~2 of \cite[Section 4]{castro2023random}-, and obtain, $\mathbb{P}$ almost surely, an $m$-mod $0$ partition $\mathcal{P}(\omega)$ of $\Delta$ and a return-time function $ R\colon  \Delta \times \Omega \to\mathbb{N}$
such that  $f^{R}_{\omega}$ is a random Gibbs-Markov map (see Definition \ref{def:random_markov_map}). In Remark \ref{measurability}, we show how we can modify the construction in Section $4.2$ in  \cite{castro2023random} to show that the random partition $\mathcal{P}(\omega)$ can be chosen measurably, and hence that $R$ is measurable. 

By \cite[Remark 4.14]{castro2023random}, there exists $\tilde{N} \in \mathbb{N}$ such that  for $\mathbb P$-almost every $\omega\in \Omega$ and $n \ge \tilde{N}$ there exists a set $E_n(\omega)\subset \Delta,$
with 
$$
m(E_n(\omega)) \le K_0 e^{-\gamma_3 n},
$$
for some $K_0,\gamma_3>0$,
such that, if we define
\begin{align*}
E_n := \bigcup_{\omega \in \Omega} \ E_n(\omega) \times  \{\omega\},
\end{align*}
then
\begin{align*}
    \{  R > n  \}\subset \{ n_{\theta_1} > n\} \cup  E_n,
\end{align*} 
where
\begin{align*}
n_{\theta_1}(x,\omega):=\min \{n \ge 0 \colon |Y^{\Delta}_m(x,\omega)|\ge \theta_1 m\,\,\forall m \ge n \}
\end{align*}
    
Because of \eqref{LargeDeviationsYoungTimes}, this concludes the proof of Theorem \ref{thm2}.
\end{proof}
\begin{remark}\label{measurability}
The measurability argument in \cite[Section 4.2]{castro2023random} can be carried almost verbatim in the general multidimensional setting. We only need to modify the following steps, where the proof in \cite{castro2023random} uses the one-dimensional structure of the circle:
\begin{itemize}
\item In \cite[Proposition 4.6]{castro2023random}, we can argue the measurability of the closure of the set  $V_n^{\delta}(x,\omega)$ as in \eqref{viennedeltaomegax}, observing that, whenever $n$ is a hyperbolic time,
the corresponding orbit segments stay away from the critical set (see Definition \ref{hyperbolic times}), so that we have
\begin{align*}
V_n^{\delta}(x,\omega) = \{y \in B_{\delta \eta^{\frac{n}{2}}}(x)  \colon f^i_{\omega}(y) \in B_{\delta \eta^{\frac{n-i}{2}}}\left(f^i_{\omega}(x)\right)\,\,\,\forall i =0,\dots,n \}.
\end{align*}
\item In \cite[Proposition 4.9]{castro2023random}, we can replace the specific dyadic partition used there with any family of countable compact sets $\{Q_j\}_{j \ge 0}$ obtained from finite covers $\{\mathcal{Q}_j\}_{j \ge 0}$ such that 
\begin{align*}
M = \bigcup_{Q \in \mathcal{Q}_j} Q,\qquad \max_{Q \in \mathcal{Q}_j}\text{diam}(Q) \to 0.
\end{align*}
The variables  $X_s(\omega) \in Q_s \cap H_n(\omega)$ used in the construction can be chosen measurably by \cite{KuratowskiRyllNardzewski1965} (for the definition of $H_n$ see \cite[pag 17]{castro2023random})
\item In the final step of the proof of Corollary $4.11$ in \cite{castro2023random}, the one-dimensional estimate \[ \# \left( \overline{w_{n,\ell}(x,\omega)} \setminus w_{n,\ell}(x,\omega) \right) \leq2 \] is replaced by \[ m\bigl(\partial w_{n,\ell}(x,\omega)\bigr)=0. \]
For the definition of $w_{n,l}(x,\omega)$, see \cite[pag 15]{castro2023random}
\end{itemize}
\end{remark}

\bibliographystyle{plain} 
\bibliography{bib}
\end{document}